\documentclass[12pt, twoside, leqno]{article}
\usepackage{amsmath,amsthm}
\usepackage{amssymb,latexsym}
\usepackage{enumerate}

\newtheorem{thm}{Theorem}[section]
\newtheorem{prop}[thm]{Proposition}
\newtheorem{cor}[thm]{Corollary}
\newtheorem{lem}[thm]{Lemma}

\newtheorem{conj}[thm]{Conjecture}

\theoremstyle{definition}

\newtheorem{rem}[thm]{Remark}

\numberwithin{equation}{section}

\usepackage{mathrsfs}
\usepackage{CJK}
\usepackage{amsmath, amscd}
\usepackage{fancyhdr}
\usepackage{enumerate}

\begin{document}

\baselineskip=17pt

\title{On Shamsuddin derivations and the isotropy groups}
\author{ Dan Yan \footnote{ The author is supported by the NSF of China (Grant No. 11871241; 11601146), the China Scholarship Council and the Construct Program of the Key Discipline in Hunan Province.}\\
MOE-LCSM,\\ School of Mathematics and Statistics,\\
 Hunan Normal University, Changsha 410081, China \\
\emph{E-mail:} yan-dan-hi@163.com \\
}
\date{}

\maketitle

\renewcommand{\thefootnote}{}

\renewcommand{\thefootnote}{\arabic{footnote}}
\setcounter{footnote}{0}

\begin{abstract} In the paper, we give an affirmative answer to the conjecture in \cite{13}. We prove that a Shamsuddin derivation $D$ is simple if and only if $\operatorname{Aut}(K[x,y_1,\allowbreak\ldots,y_n])_D=\{id\}$. In addition, we calculate the isotropy groups of the Shamsuddin derivations $D=\partial_x+\sum_{j=1}^r(a(x)y_j+b_j(x))\partial_j$ of $K[x,y_1,\ldots,y_r]$. We also prove that $\operatorname{Im}D$ is a Mathieu-Zhao subspace if and only if $a(x)\in K$.
\end{abstract}
 \vskip 2.5mm

\section{Introduction}

Throughout this paper, we will write $\mathbb{N}$ for the non-negative integers,
$K$ for any field with
characteristic zero and $R:=K[x,y_1,\ldots,y_n]$ for the
polynomial algebra over $K$ in $n+1$ indeterminates $x,y_1,\ldots,y_n$.
$\partial_x,~\partial_i$ will denote the derivations $\frac{\partial}{\partial x}$, $\frac{\partial}{\partial y_i}$ of $R$ for all $1\leq i\leq n$, respectively. More generally, if $s,r_1,\ldots, r_s\geq 1$ are integers and $\{x\}\bigcup\{y_{i,j}:i=1,\ldots,s,~j=1,\ldots,r_i\}$ are indeterminates over $K$, $\partial_{i,j}$ will denote the derivation $\frac{\partial}{\partial y_{i,j}}$ of $K[x,\bigcup_{i=1}^s\{y_{i,1},\ldots,y_{i,r_i}\}]$. We abbreviate $\frac{\partial g_t}{\partial y_j}$ as $g_{ty_j}$. For element $f$ of $K[x]$, we shall often use $f'$ instead of $f_x$.

A $K$-derivation $D:R\rightarrow R$ of $R$ is a $
K$-linear map such that
$$D(ab)=D(a)b+aD(b)$$
for any $a,b\in R$ and $D(c)=0$ for any $c\in K$. The set of all $K$-derivations of $R$
is denoted by $\operatorname{Der}_ K(R)$. An ideal $I$ of $R$ is called $D$-$stable$ if $D(I)\subset I$. $R$ is called
$D$-$simple$ if it has no proper nonzero $D$-$stable$ ideal. The $K$-derivation $D$ is called $simple$ if $R$ has no $D$-$stable$
ideals other than $0$ and $R$. For some examples of simple
derivations, see \cite{4}, \cite{9}, \cite{5}, \cite{7}, \cite{2}.

Let
$\operatorname{Aut}(R)$ act on $\operatorname{Der}_ K( R)$
by:
$$(\rho,D)\rightarrow\rho^{-1}\circ D\circ\rho=\rho^{-1} D\rho.$$
The isotropy subgroup is defined to be:
$$\operatorname{Aut}(R)_D:=\{\rho\in\operatorname{Aut}(R)|\rho D=D\rho\}.$$

A derivation $D$ of $R$ is said to be a Shamsuddin derivation if
$D=\partial_x+\sum_{i=1}^n(a_iy_i+b_i)\partial_i$ with $a_i,
b_i\in  K[x]$ for all $1\leq i\leq n$. Observe that if $D$ is such a Shamsuddin derivation of $R$, then grouping the terms that have the same $a_i$ and renaming the indeterminates $y_i$ and the polynomials $a_i, b_i$ if necessary, we can write $D$ in the following form:
$$D=\partial_x+\sum_{i=1}^s\sum_{j=1}^{r_i}(a_iy_{i,j}+b_{i,j})\partial_{i,j}$$
with $a_i, b_{i,j}\in K[x]$ for every $i$ and every $(i,j)$, $a_i\neq a_l$ for $i\neq l$. A derivation $D$ of $R$ is said to be locally finite if, for each $a\in R$, the $K$-subspace spanned by $D^i(a)$ $(i\geq 1)$ is finite dimension over $K$.

The Mathieu-Zhao subspace was introduced by Zhao in \cite{12} and \cite{14}, which is a natural generalization of ideals. We give the definition here for the polynomial rings. A $K$-subspace $M$ of $R$ is said to be a Mathieu-Zhao subspace if for any $a, b\in R$ with $a^m\in M$ for all $m\geq 1$, we have $ba^m\in M$ when $m>>0$.

In \cite{13}, L.N.Bertoncello and
D.Levcovitz have proved that the isotropy group of simple Shamsuddin derivations is trivial. They also proposed the following conjecture.

\begin{conj}
 If $D$ is a Shamsuddin derivation
of $K[x_1,x_2,\ldots,x_n]$, then $D$ is simple if, and only if $\operatorname{Aut}( K[x_1,x_2,\ldots,x_n])_D\allowbreak =\{id\}$.
\end{conj}

In our paper, we give an affirmative answer to the conjecture. In addition, we calculate the isotropy groups of the Shamsuddin derivations
$D=\partial_x+\sum_{j=1}^r(a(x)y_j+b_j(x))\partial_j$ of $K[x,y_1,\ldots,y_r]$. In section 3, we also prove that $\operatorname{Im}D$ is a Mathieu-Zhao subspace if and only if $a(x)\in K$.

\section{Affirmative answer to the conjecture in \cite{13}}

\begin{lem} \label{lem2.1}
Let $D=\partial_x+\sum_{j=1}^rb_j(x)\partial_j$ be a derivation of $K[x,y_1,\ldots,y_r]$. Then $\operatorname{Aut}(K[x,y_1,\ldots,y_r])_D=\{(f,g_1,\ldots,g_r)\}$, where
$f=x+p(y_1-h_1(x),\ldots,y_r-h_r(x))$, $g_t=\sum_{k=0}^{m_t}\frac{1}{k+1}b_{t,k}f^{k+1}+q_t(y_1-h_1(x),\ldots,y_r-h_r(x))$,
where $b_t(x)=\sum_{k=0}^{m_t}b_{t,k}x^k$, $h_t(x)=\int b_t(x)dx$ for all $1\leq t\leq r$, $p\in K[y_1-h_1(x),\ldots,y_r-h_r(x)]$, $(f,q_1,\ldots,q_r)$ is any polynomial automorphism of $K[x,y_1,\ldots,y_r]$.
\end{lem}
\begin{proof}
Let $\rho\in \operatorname{Aut}(K[x,y_1,\ldots,y_r])_D$ with $\rho(x)=f(x,y_1,\ldots,y_r)$, $\rho(y_t)=g_t(x,y_1,\allowbreak\ldots,y_r)$ for all $1\leq t\leq r$. Then we have the following equations:
\begin{eqnarray}
  D(\rho(x))=\rho(D(x))\label{eq2.1}\\
  D(\rho(y_t))=\rho(D(y_t))\label{eq2.2}
\end{eqnarray}
for all $1\leq t\leq r$. That is,
\begin{eqnarray}
  f_x+\sum_{j=1}^rb_j(x)f_{y_j}=1\label{eq2.3}\\
  g_{tx}+\sum_{j=1}^rb_j(x)g_{ty_j}=b_t(f)\label{eq2.4}
\end{eqnarray}
for all $1\leq t\leq r$. Let $\bar{x}=x$, $\bar{y}_j=y_j-h_j(x)$, where $h_j(x)=\int b_j(x)dx$ for all $1\leq j\leq r$. Then it follows from equation \eqref{eq2.3} that $f_{\bar{x}}=1$. Thus, we have $f=\bar{x}+p(\bar{y}_1,\ldots,\bar{y}_r)$ for some $p\in K[\bar{y}_1,\ldots,\bar{y}_r]$. That is, $f=x+p(y_1-h_1(x),\ldots,y_r-h_r(x))$. Since $b_t(x)=\sum_{k=0}^{m_t}b_{t,k}x^k$, it follows from equation \eqref{eq2.4} that $g_{t\bar{x}}=\sum_{k=0}^{m_t}b_{t,k}f^k$. Thus, we have
$$g_t=\sum_{k=0}^{m_t}\frac{1}{k+1}b_{t,k}f^{k+1}+q_t(\bar{y}_1,\ldots,\bar{y}_r)$$
for some $q_t(\bar{y}_1,\ldots,\bar{y}_r)\in K[\bar{y}_1,\ldots,\bar{y}_r]$ and for all $1\leq t\leq r$. Since $\rho \in \operatorname{Aut}(K[x,y_1,\allowbreak\ldots,y_r])$, we have that $(f,q_1,\ldots,q_r)$ is a polynomial automorphism of $K[x,y_1,\ldots,\allowbreak y_r]$, which completes the proof.
\end{proof}

\begin{rem}\label{rem2.2}
Let $D=\partial_x+\sum_{j=1}^r(a(x)y_j+b_j(x))\partial_j$ be a derivation of $K[x,y_1,\allowbreak\ldots,y_r]$. If there exists $j_0\in \{1,2,\ldots,r\}$ such that $b_{j_0}(x)=0$, then it is easy to verify that $(x,y_1,\ldots,\tilde{c}y_{j_0},\ldots,y_r)\in \operatorname{Aut}(K[x,y_1,\allowbreak\ldots,y_r])_D$ for any $\tilde{c}\in K^*$. Thus, $\operatorname{Aut}(K[x,y_1,\allowbreak\ldots,y_r])_D \neq \{id\}$. If $a(x)=0$, then let $(f,q_1,\ldots,q_r)=(x,\tilde{c}_1(y_1-\sum_{k=0}^{m_1}\frac{1}{k+1}b_{1,k}x^{k+1}),\allowbreak\ldots,\tilde{c}_r(y_r-\sum_{k=0}^{m_r}\frac{1}{k+1}b_{r,k}x^{k+1}))$ in Lemma \ref{lem2.1} with $\tilde{c}_1,\ldots,\tilde{c}_r\in K^*$. We have $(x,\sum_{k=0}^{m_1}\frac{1}{k+1}b_{1,k}x^{k+1}+\tilde{c}_1(y_1-\sum_{k=0}^{m_1}\frac{1}{k+1}b_{1,k}x^{k+1}),\ldots,\sum_{k=0}^{m_r}\frac{1}{k+1}b_{r,k}x^{k+1}\allowbreak+\tilde{c}_r(y_r-\sum_{k=0}^{m_r}\allowbreak\frac{1}{k+1}b_{r,k}x^{k+1}))\in \operatorname{Aut}(K[x,y_1,\ldots,y_r])_D$ for any $\tilde{c}_1,\ldots,\tilde{c}_r\in K^*$. Hence $\operatorname{Aut}(K[x,y_1,\allowbreak\ldots,y_r])_D\neq \{id\}$ in Lemma \ref{lem2.1}.
\end{rem}

\begin{thm}\label{thm2.3}
Let $D=\partial_x+\sum_{j=1}^r(a(x)y_j+b_j(x))\partial_j$ be a Shamsuddin derivation of $K[x,y_1,\ldots,y_r]$. Then the following two statements are equivalent:

$(1)$ $D$ is a simple derivation;

$(2)$ $\operatorname{Aut}(K[x,y_1,\ldots,y_r])_D= \{id\}$.
\end{thm}
\begin{proof}
$(1)\Rightarrow (2)$ It follows from Theorem 3.2 in \cite{13}.

$(2)\Rightarrow (1)$ It follows from Theorem 3.2 in \cite{6} that $D$ is a simple derivation if and only if $z'=a(x)z+\sum_{j=1}^rk_jb_j(x)$ does not have any solution in $K[x]$ for every $(k_1,\ldots,k_r)\in K^r\setminus\{(0,\ldots,0)\}$. Thus, it suffices to prove that if $z'=a(x)z+\sum_{j=1}^rk_jb_j(x)$ has a solution in $K[x]$ for some $(k_1,\ldots,k_r)\in K^r\setminus\{(0,\ldots,0)\}$, then $\operatorname{Aut}(K[x,y_1,\ldots,y_r])_D\neq \{id\}$. It follows from Lemma \ref{lem2.1} and Remark \ref{rem2.2} that we can assume that $a(x)b_1(x)\cdots b_r(x)\neq 0$. Let $\rho \in \operatorname{Aut}(K[x,y_1,\ldots,y_r])_D$ and $\rho(x)=f(x,y_1,\ldots,y_r)$, $\rho(y_t)=g_t(x,y_1,\ldots,y_r)$ for \allowbreak all $1\leq t\leq r$; $f(x,y_1,\ldots,y_r)=\sum_{|\alpha|=d}f_{\alpha}(x)y_1^{\alpha_1}\cdots y_r^{\alpha_r}+\sum_{|\alpha|<d}f_{\alpha}(x)y_1^{\alpha_1}\cdots y_r^{\alpha_r}$, $g_t(x,y_1,\ldots,y_r)=\sum_{|\beta_t|=n_t}g_{t,\beta_t}(x)y_1^{\beta_{t1}}\cdots y_r^{\beta_{tr}}+\sum_{|\beta_t|<n_t}g_{t,\beta_t}(x)y_1^{\beta_{t1}}\cdots y_r^{\beta_{tr}}$ with $f_{\alpha}(x)\allowbreak\neq 0$, $g_{t,\beta_t}(x)\neq 0$ for some $|\alpha|=\alpha_1+\cdots+\alpha_r$, $|\beta_t|=\beta_{t1}+\cdots+\beta_{tr}$ and for all $1\leq t\leq r$. Then we have the following equations:
\begin{eqnarray}
  D(\rho(x))=\rho(D(x))\label{eq2.5}\\
  D(\rho(y_t))=\rho(D(y_t))\label{eq2.6}
\end{eqnarray}
for all $1\leq t\leq r$. It follows from equation \eqref{eq2.5} that
\begin{equation}\label{eq2.7}
\begin{split}
\sum_{|\alpha|=d}f'_{\alpha}(x)y_1^{\alpha_1}\cdots y_r^{\alpha_r}+\sum_{|\alpha|<d}f'_{\alpha}(x)y_1^{\alpha_1}\cdots y_r^{\alpha_r}+\sum_{j=1}^r(a(x)y_j+b_j(x))\times \\
(\sum_{|\alpha|=d}\alpha_jf_{\alpha}(x)y_1^{\alpha_1}\cdots y_j^{\alpha_j-1}\cdots y_r^{\alpha_r}+\sum_{|\alpha|<d}\alpha_jf_{\alpha}(x)y_1^{\alpha_1}\cdots y_j^{\alpha_j-1}\cdots y_r^{\alpha_r})=1
\end{split}
\end{equation}
We view the polynomials in $K[x][y_1,\ldots,y_r]$ with coefficients in $K[x]$ when we compare the coefficients of $y_1^{\alpha_1}\cdots y_r^{\alpha_r}$. If $d \geq 1$, then we have
\begin{eqnarray}\label{eq2.8}
 f'_{\alpha}(x)=-(\alpha_1+\cdots+\alpha_r)a(x)f_{\alpha}(x)
\end{eqnarray}
by comparing the coefficients of $y_1^{\alpha_1}\cdots y_r^{\alpha_r}$ with $|\alpha|=d$ of equation \eqref{eq2.7}. Thus, we have $f_{\alpha}(x)=0$ for all $|\alpha|=d$ by comparing the degree of $x$ of equation \eqref{eq2.8}, which is a contradiction. Therefore, we have $d=0$. That is, $f(x,y_1,\ldots,y_r)=f_0(x)$. It follows from equation \eqref{eq2.3} that $f'_0(x)=1$. Thus, we have $f_0(x)=x+c$ for some $c\in K$. That is, $f(x,y_1,\ldots,y_r)=x+c$. It follows from equation \eqref{eq2.6} that
\begin{equation}\label{eq2.9}
\begin{split}
\sum_{|\beta_t|=n_t}g'_{t,\beta_t}(x)y_1^{\beta_{t1}}\cdots y_r^{\beta_{tr}}+\sum_{|\beta_t|<n_t}g'_{t,\beta_t}(x)y_1^{\beta_{t1}}\cdots y_r^{\beta_{tr}}+\sum_{j=1}^r(a(x)y_j+b_j(x))\cdot \\
(\sum_{|\beta_t|=n_t}\beta_{tj}\cdot g_{t,\beta_t}(x)y_1^{\beta_{t1}}\cdots y_j^{\beta_{tj}-1}\cdots y_r^{\beta_{tr}}+\sum_{|\beta_t|<n_t}\beta_{tj}\cdot g_{t,\beta_t}(x)y_1^{\beta_{t1}}\cdots y_j^{\beta_{tj}-1}\cdots y_r^{\beta_{tr}})\\=a(x+c)(\sum_{|\beta_t|=n_t}g_{t,\beta_t}(x)y_1^{\beta_{t1}}\cdots y_r^{\beta_{tr}}+\sum_{|\beta_t|<n_t}g_{t,\beta_t}(x)y_1^{\beta_{t1}}\cdots y_r^{\beta_{tr}})+b_t(x+c)
\end{split}
\end{equation}
If $n_t \geq 1$, then we have
\begin{eqnarray}\label{eq2.10}
 g'_{t,\beta_t}(x)=[a(x+c)-(\beta_{t1}+\cdots+\beta_{tr})a(x)]g_{t,\beta_t}(x)
\end{eqnarray}
by comparing the coefficients of $y_1^{\beta_{t1}}\cdots y_r^{\beta_{tr}}$ with $|\beta_t|=n_t$ of equation \eqref{eq2.9} for all $1\leq t\leq r$. Thus, we have
$g_{t,\beta_t}(x)\in K^*$ for all $|\beta_t|=n_t$ and
\begin{eqnarray}\label{eq2.11}
 a(x+c)=(\beta_{t1}+\cdots+\beta_{tr})a(x)
\end{eqnarray}
by comparing the degree of $x$ of equation \eqref{eq2.10}. Thus, we have $|\beta_t|=1$ by comparing the highest degree of $x$ of equation \eqref{eq2.11} for all $1\leq t\leq r$. Therefore, there exists $j_t\in \{1,2,\ldots,r\}$ such that $\beta_{tj_t}=1$ and $\beta_{tj}=0$ for all $j\neq j_t$ and for all $1\leq t\leq r$. Therefore, we have
\begin{eqnarray}\label{eq2.12}
 g_t=\sum_{j=1}^rc_{tj}y_j+g_{t,0}(x)
\end{eqnarray}
for some $c_{tj}\in K$, $\det(c_{tj})_{r\times r}\neq 0$ and for all $1\leq t,j\leq r$.

(1) If $c\neq 0$, then $\operatorname{Aut}(K[x,y_1,\ldots,y_r])_D\neq\{id\}$ because $(x+c,g_1,\ldots,g_r)\in \operatorname{Aut}(K[x,y_1,\ldots,y_r])_D$ for some $g_1,\ldots,g_t\in K[x,y_1,\ldots,y_r]$.

(2) If $c=0$, then it follows from equation \eqref{eq2.9} that
\begin{eqnarray}\label{eq2.13}
 g'_{t,0}(x)=a(x)g_{t,0}(x)+b_t(x)-\sum_{j=1}^rc_{tj}b_j(x)
\end{eqnarray}
for all $1\leq t\leq r$.

Suppose that $Q(x)\in K[x]$ is a solution of $z'=a(x)z+\sum_{j=1}^rk_jb_j(x)$ for some $(k_1,\ldots,k_r)\in K^r\setminus \{(0,\ldots,0)\}$. Without loss of generality, we can assume that $k_1=1$. Since equation \eqref{eq2.13} is satisfied by $g_{t,0}(x)$ for all $1\leq t\leq r$, we can choose a $t$ such that $t=1$. Let $1-c_{11}=e$, $c_{1j}=-e\cdot k_j$ for all $2\leq j\leq r$. Then $g_{1,0}(x)=e\cdot Q(x)$ satisfies equation \eqref{eq2.13} for any $e\in K^*$. Thus, $g_1(x,y_1,\ldots,y_r)=(1-e)y_1-e\cdot \sum_{j=2}^rk_jy_j+e\cdot Q(x)$ for any $e\in K^*$. Therefore, we have $(x,(1-e)y_1-e\cdot \sum_{j=2}^rk_jy_j+e\cdot Q(x),y_2,\ldots,y_r)\in \operatorname{Aut}(K[x,y_1,\ldots,y_r])_D$ for any $e\in K^*$, $e\neq 1$. That is, $\operatorname{Aut}(K[x,y_1,\ldots,y_r])_D\neq\{id\}$, which completes the proof.
\end{proof}

\begin{cor}\label{cor2.4}
Let $D=\partial_x+\sum_{j=1}^r(a(x)y_j+b_j(x))\partial_j$ be a Shamsuddin derivation of $K[x,y_1,\ldots,y_r]$ with $a(x)\neq 0$. Then $\operatorname{Aut}(K[x,y_1,\ldots,y_r])_D=\{(x+c,\sum_{j=1}^rc_{1j}y_j+g_{1,0}(x),\ldots,\sum_{j=1}^rc_{rj}y_j+g_{r,0}(x))|\det(c_{tj})_{r\times r}\neq 0~ and~ g'_{t,0}(x)=a(x)\cdot \allowbreak g_{t,0}(x)+b_t(x+c)-\sum_{j=1}^rc_{tj}b_j(x)~ for~ all~ 1\leq t\leq r\}$.
\end{cor}
\begin{proof}
The conclusion follows from the proof of Theorem \ref{thm2.3}.
\end{proof}

\begin{prop}\label{prop2.5}
Let $D=\partial_x+\sum_{j=1}^r(a(x)y_j+b_j(x))\partial_j$ be a Shamsuddin derivation of $K[x,y_1,\ldots,y_r]$ with $a(x)\neq 0$. Then we have the following statements:

$(1)$ If $\operatorname{deg}a(x)\geq 1$, then $\operatorname{Aut}(K[x,y_1,\ldots,y_r])_D=\{(x,\sum_{j=1}^rc_{1j}y_j+g_{1,0}(x),\ldots,\allowbreak \sum_{j=1}^rc_{rj}y_j+g_{r,0}(x))|\det(c_{tj})_{r\times r}\neq 0~ and~ g'_{t,0}(x)=a(x)g_{t,0}(x)+b_t(x)-\sum_{j=1}^rc_{tj}b_j(x)\allowbreak~ for~ all~ 1\leq t\leq r\}$.

$(2)$ If $a(x)\in K^*$, then $\operatorname{Aut}(K[x,y_1,\ldots,y_r])_D=\{(x+c,\sum_{j=1}^rc_{1j}y_j+g_{1,0}(x),\ldots,\allowbreak\sum_{j=1}^rc_{rj}y_j+g_{r,0}(x))|\det(c_{tj})_{r\times r}\neq 0~ and~ g_{t,0}(x)=\hat{c}e^{ax}+(\int f_t(x)e^{-ax}dx)e^{ax}~ for~ all~\allowbreak 1\leq t\leq r\}$, where $a:=a(x)$, $\hat{c}\in K$ and $f_t(x)=b_t(x+c)-\sum_{j=1}^rc_{tj}b_j(x)$ for $1\leq t\leq r$.
\end{prop}
\begin{proof}
It follows from Corollary \ref{cor2.4} that $\operatorname{Aut}(K[x,y_1,\ldots,y_r])_D=\{(x+c,\sum_{j=1}^rc_{1j}y_j\allowbreak+g_{1,0}(x),\ldots,\sum_{j=1}^rc_{rj}y_j+g_{r,0}(x))|\det(c_{tj})_{r\times r}\neq 0~ and~ g'_{t,0}(x)=a(x)g_{t,0}(x)+b_t(x+c)-\sum_{j=1}^rc_{tj}b_j(x)~ for~ all~ 1\leq t\leq r\}$.

$(1)$ It follows from equation \eqref{eq2.11} and the arguments in Theorem \ref{thm2.3} that
\begin{eqnarray}\label{eq2.14}
a(x+c)=a(x).
\end{eqnarray}
If $\operatorname{deg}a(x)\geq 1$, then we have $c=0$ by comparing the coefficients of the monomial $x^{\operatorname{deg}a(x)-1}$ of equation \eqref{eq2.14}. Then the conclusion follows.

$(2)$ If $a(x)\in K^*$, then let $a:=a(x)$, it follows from equation \eqref{eq2.9} and the arguments in Theorem \ref{thm2.3} that
\begin{eqnarray}\label{eq2.15}
 g'_{t,0}(x)=ag_{t,0}(x)+b_t(x+c)-\sum_{j=1}^rc_{tj}b_j(x)
\end{eqnarray}
for all $1\leq t\leq r$.
Let $f_t(x)=b_t(x+c)-\sum_{j=1}^rc_{tj}b_j(x)$ for $1\leq t\leq r$. Then it's easy to compute that $g_{t,0}(x)=\hat{c}e^{ax}+(\int f_t(x)e^{-ax}dx)e^{ax}$ for any $\hat{c}\in K$ for all $1\leq t\leq r$. Then the conclusion follows.
\end{proof}

\begin{rem}\label{rem2.6}
If $a(x)\in K^*$ in Proposition \ref{prop2.5}, then it follows from Theorem 3.2 (a) in \cite{6} that $D$ is not simple. Since $f_1(x),\ldots,f_r(x)$ in Proposition \ref{prop2.5} are polynomials, we can see from Proposition \ref{prop2.5} (2) that there are polynomials $g_{1,0}(x),\ldots,g_{r,0}(x)$ such that $(x,\sum_{j=1}^rc_{1j}y_j+g_{1,0}(x),\ldots,\sum_{j=1}^rc_{rj}y_j+g_{r,0}(x))\in \operatorname{Aut}(K[x,y_1,\ldots,y_r])_D$ for any $\det(c_{tj})_{r\times r}\neq 0$.
\end{rem}

\begin{thm}\label{thm2.7}
Let $D=\partial_x+\sum_{i=1}^s\sum_{j=1}^{r_i}(a_i(x)y_{i,j}+b_{i,j}(x))\partial_{i,j}$ be a Shamsuddin derivation of $K[x,y_{1,1},\ldots,y_{1,r_1},\ldots,y_{s,1},\ldots,y_{s,r_s}]$. If $A_1=(x, g_{i_0,1},\ldots,g_{i_0,r_{i_0}})\in \operatorname{Aut}(K[x,y_{i_0,1},\ldots,y_{i_0,r_{i_0}}])_{D_{i_0}}$ for some $i_0\in \{1,2,\ldots,s\}$, where $D_{i_0}=\partial_x+\sum_{j=1}^{r_{i_0}}\allowbreak (a_{i_0}(x)y_{i_0,j}+b_{i_0,j}(x))\partial_{i_0,j}$, then $A_2=(x,y_{1,1},\ldots,y_{1,r_1},\allowbreak\ldots,y_{i_0-1,1},\ldots,y_{i_0-1,r_{i_0-1}},\allowbreak g_{i_0,1},\ldots,g_{i_0,r_{i_0}},y_{i_0+1,1},\ldots,y_{i_0+1,r_{i_0+1}},\ldots,y_{s,1},\ldots,y_{s,r_s})\in \operatorname{Aut}(K[x,y_1,\ldots,y_n])_D$.
\end{thm}
\begin{proof}
Without loss of generality, we can assume that $i_0=1$. Clearly, $A_2$ is a polynomial automorphism of $K[x,y_1,\ldots,y_n]$. Let $\rho=A_2$. That is, $\rho(x)=x$, $\rho(y_{1,j})=g_{1,j}$ and $\rho(y_{i,j})=y_{i,j}$ for all $1\leq j\leq r_i$, $2\leq i\leq s$. It suffices to prove that $D(\rho(x))=\rho(D(x))$ and $D(\rho(y_{i,j}))=\rho(D(y_{i,j}))$ for all $1\leq j\leq r_i$, $1\leq i\leq s$. Since $A_1\in \operatorname{Aut}(K[x,y_{1,1},\ldots,y_{1,r_1}])_{D_1}$, we have $D_1(\rho(x))=\rho(D_1(x))$ and $D_1(\rho(y_{1,j}))=\rho(D_1(y_{1,j}))$ for all $1\leq j\leq r_1$. That is,
\begin{eqnarray}\label{eq2.16}
 (g_{1,j})_x+\sum_{j=1}^{r_1}(a_1(x)y_{1,j}+b_{1,j}(x))(g_{1,j})_{y_{1,j}}=a_1(f)g_{1,j}+b_{1,j}(f)
\end{eqnarray}
for all $1\leq j\leq r_1$. Since $D(\rho(x))=1$, $\rho(D(x))=1$, $D(\rho(y_{1,j}))=(g_{1,j})_x+\sum_{j=1}^{r_1}(a_1(x)y_{1,j}+b_{1,j}(x))(g_{1,j})_{y_{1,j}}$
and $\rho(D(y_{1,j}))=a_1(f)g_{1,j}+b_{1,j}(f)$, $D(\rho(y_{i,j}))=a_iy_{i,j}+b_{i,j}=\rho(D(y_{i,j}))$
for all $1\leq j\leq r_i$, $2\leq i\leq s$. It follows from equation \eqref{eq2.16} that $D(\rho(x))=\rho(D(x))$ and $D(\rho(y_{i,j}))=\rho(D(y_{i,j}))$ for all $1\leq j\leq r_i$, $1\leq i\leq s$. That is, $A_2\in \operatorname{Aut}(K[x,y_1,\ldots,y_n])_D$, which completes the proof.
\end{proof}

\begin{thm}\label{thm2.8}
Let $D=\partial_x+\sum_{i=1}^s\sum_{j=1}^{r_i}(a_i(x)y_{i,j}+b_{i,j}(x))\partial_{i,j}$ be a Shamsuddin derivation of $K[x,y_1,\ldots,y_n]$. Then the following two statements are equivalent:

(1) $D$ is a simple derivation;

(2) $\operatorname{Aut}(K[x,y_1,\ldots,y_n])_D=\{id\}$.
\end{thm}
\begin{proof}
$(1)\Rightarrow (2)$ It follows from Theorem 3.2 in \cite{13}.

$(2)\Rightarrow (1)$ It suffices to show that if $D$ is not simple, then $\operatorname{Aut}(K[x,y_1,\ldots,y_n])_D\allowbreak\neq\{id\}$. Since $D$ is not simple, it follows from Theorem 3.1 in \cite{6} that $D_{i_0}$ is not simple for some $i_0\in \{1,2,\ldots,s\}$, where $D_{i_0}=\partial_x+\sum_{j=1}^{r_{i_0}}(a_{i_0}(x)y_{i_0,j}+b_{i_0,j}(x))\partial_{i_0,j}$.
Without loss of generality, we can assume that $i_0=1$.

If $\operatorname{deg}a_1(x)\geq 1$, then it follows from Proposition \ref{prop2.5} (1) and Theorem \ref{thm2.3} that there exists $A_1\in \operatorname{Aut}(K[x,y_{1,1},\ldots,y_{1,r_1}])$ with $\rho(x)=x$ and $A_1\neq id$ such that $A_1\in \operatorname{Aut}(K[x,y_{1,1},\ldots,y_{1,r_1}])_{D_1}$.

If $a_1(x)=0$, then it follows from Remark \ref{rem2.2} that there exists $A_1\in \operatorname{Aut}(K[x,\allowbreak y_{1,1},\ldots,y_{1,r_1}])$ with $\rho(x)=x$ and $A_1\neq id$ such that $A_1\in \operatorname{Aut}(K[x,y_{1,1},\ldots,y_{1,r_1}])_{D_1}$.

If $a_1(x)\in K^*$, then it follows from Proposition \ref{prop2.5} (2) and Remark \ref{rem2.6} that there exists $A_1\in \operatorname{Aut}(K[x,y_{1,1},\ldots,y_{1,r_1}])$ with $\rho(x)=x$ and $A_1\neq id$ such that $A_1\in \operatorname{Aut}(K[x,y_{1,1},\ldots,y_{1,r_1}])_{D_1}$.

It follows from Theorem \ref{thm2.7}  that $A_2=(A_1,y_{2,1},\ldots,y_{2,r_2},\ldots, y_{s,1},\ldots,y_{s,r_s})\in \operatorname{Aut}(K[x,y_1,\ldots,y_n])_D$. Since $A_1\neq id$, we have $A_2\neq id$. Then the conclusion follows.
\end{proof}

\begin{rem}\label{rem2.9}
In \cite{11}, the authors have proved Theorem \ref{thm2.8} if $n=1$ and $a_1(x)\neq 0$. We have proved Theorem \ref{thm2.8} if $n=1$ in \cite{3}.
\end{rem}

\section{Images of Shamsuddin derivations}

\begin{lem}\label{lem3.1}
Let $D=\partial_x+\sum_{j=1}^n(a_j(x)y_j+b_j(x,y_1,\ldots,y_{j-1}))\partial_j$ be a derivation of $K[x,y_1,\ldots,y_n]$. Then $D$ is locally finite if and only if $a_j(x)\in K$ for all $1\leq j\leq n$.
\end{lem}
\begin{proof}
$``\Leftarrow"$ The conclusion follows from Example 9.3.2 in \cite{10}.

$``\Rightarrow"$ It suffices to prove that if there exists $i_0\in \{1,2,\ldots,n\}$ such that $\deg a_{i_0}(x)\geq 1$, then $D$ is not locally finite. Since $D(y_{i_0})=a_{i_0}(x)y_{i_0}+b_{i_0}(x,y_1,\ldots,\allowbreak y_{i_0-1})$, we have
$$D^2(y_{i_0})=(a^2_{i_0}(x)+a'_{i_0}(x))y_{i_0}+\operatorname{polynomial}~\operatorname{in}~ K[x,y_1,\ldots,y_{i_0-1}].$$
Suppose that $D^{k-1}(y_{i_0})=a^{k-1}_{i_0}(x)y_{i_0}+P_{k-1}(x)y_{i_0}+\operatorname{polynomial}~\operatorname{in}~ K[x,y_1,\ldots,\allowbreak y_{i_0-1}]$, where $P_{k-1}(x)$ is a polynomial in $K[x]$ and $\deg P_{k-1}(x)<(k-1)\deg a_{i_0}(x)$. Then we have
$$D^k(y_{i_0})=D(D^{k-1}(y_{i_0}))=a^k_{i_0}(x)y_{i_0}+P_k(x)y_{i_0}+\operatorname{polynomial}~\operatorname{in}~ K[x,y_1,\ldots, y_{i_0-1}],$$
where  $P_k(x)$ is a polynomial in $K[x]$ and $\deg P_k(x)<k\deg a_{i_0}(x)$. Let $S=\{D^i(y_{i_0})|i\geq 1\}$. Then the vector space generated by $S$ is infinite dimension over $K$. Thus, $D$ is not locally finite, which completes the proof.
\end{proof}

\begin{thm}\label{thm3.2}
Let $D=\partial_x+\sum_{j=1}^n(a_j(x)y_j+b_j(x))\partial_j$ be a derivation of $K[x,y_1,\allowbreak\ldots,y_n]$. If the equation $\sum_{i=1}^n\gamma_ia_i(x)=0$ has no non-zero solutions in $\mathbb{N}^n$ for $\gamma=(\gamma_1,\gamma_2,\ldots,\gamma_n)$ and there exists $i_0\in\{1,2,\allowbreak\ldots,n\}$ such that $\deg a_{i_0}(x)\geq 1$, then $\operatorname{Im}D$ is not a Mathieu-Zhao subspace.
\end{thm}
\begin{proof}
Since $1\in \operatorname{Im}D$, we have $\operatorname{Im}D=K[x,y_1,\ldots,y_n]$ if $\operatorname{Im}D$ is a Mathieu-Zhao subspace. Without loss of generality, we assume that $i_0=1$. We claim that $y_1\notin \operatorname{Im}D$. Suppose that $y_1\in \operatorname{Im}D$. Then there exists $f(x,y_1,\ldots,y_n)\in K[x,y_1,\ldots,y_n]$ such that
\begin{eqnarray}
  D(f(x,y_1,\ldots,y_n))=y_1.\label{eq3.1}
\end{eqnarray}
Let $f=f^{(d)}+f^{(d-1)}+\cdots+f^{(1)}+f^{(0)}$ with $f^{(d)}\neq 0$, where $f^{(i)}$ is the homogeneous part of degree $i$ with respect to $y_1,\ldots,y_n$ of $f$. It follows from equation \eqref{eq3.1} that
\begin{eqnarray}
  f_x+\sum_{j=1}^n(a_j(x)y_j+b_j(x))f_{y_j}=y_1.\label{eq3.2}
\end{eqnarray}
If $d\geq 2$, then we have
\begin{eqnarray}
  f_x^{(d)}+\sum_{j=1}^na_j(x)y_jf_{y_j}^{(d)}=0.\label{eq3.3}
\end{eqnarray}
by comparing the homogeneous part of degree $d$ with respect to $y_1,\ldots,y_n$ of equation \eqref{eq3.2}. Let $f^{(d)}=\sum_{l_1+\cdots+l_n=d}a_{l_1\cdots l_n}y_1^{l_1}\cdots y_n^{l_n}$ with $a_{l_1\cdots l_n}=a_{l_1\cdots l_n}(x)\in K[x]$. Then equation \eqref{eq3.3} has the following form
\begin{eqnarray}
 \sum_{l_1+\cdots+l_n=d}a_{l_1\cdots l_n}'(x)y_1^{l_1}\cdots y_n^{l_n}+\sum_{j=1}^nl_ja_j(x)\cdot\sum_{l_1+\cdots+l_n=d}a_{l_1\cdots l_n}(x)y_1^{l_1}\cdots y_n^{l_n}=0.\label{eq3.4}
\end{eqnarray}
We view the polynomials in $K[x][y_1,\ldots,y_n]$ with coefficients in $K[x]$ when we comparing the coefficients of $y_1^{l_1}\cdots y_n^{l_n}$. Thus, we have
\begin{eqnarray}
 a_{l_1\cdots l_n}'(x)+(\sum_{j=1}^nl_ja_j(x))a_{l_1\cdots l_n}(x)=0\label{eq3.5}
\end{eqnarray}
by comparing the coefficients of $y_1^{l_1}\cdots y_n^{l_n}$ of equation \eqref{eq3.4} with $l_1+\cdots+l_n=d$. Since $(l_1,l_2,\ldots,l_n)\in \mathbb{N}^n$ and $(l_1,l_2,\ldots,l_n)\neq (0,0,\ldots,0)$, we have $\sum_{j=1}^nl_ja_j(x)\neq 0$. Thus, we have $a_{l_1\cdots l_n}(x)=0$ for all $l_1+\cdots+l_n=d$ by comparing the degree of $x$ of equation \eqref{eq3.5}. That is, $f^{(d)}=0$, which is a contradiction. Thus, we have $d\leq 1$.

If $d=1$, then $f^{(1)}=c_1(x)y_1+\cdots+c_n(x)y_n$ with $c_j(x)\in K[x]$ for $1\leq j\leq n$. Thus, equation \eqref{eq3.2} has the following form
\begin{eqnarray}
 f_x^{(0)}+\sum_{j=1}^nc_j'(x)y_j+\sum_{j=1}^n(a_j(x)y_j+b_j(x))c_j(x)=y_1.\label{eq3.6}
\end{eqnarray}
Hence, we have
\begin{eqnarray}
 c_1'(x)+a_1(x)c_1(x)=1\label{eq3.7}
\end{eqnarray}
by comparing the coefficients of $y_1$ of equation \eqref{eq3.6}. Since $\deg a_1(x)\geq 1$, we have $c_1(x)=0$ by comparing the degree of $x$ of equation \eqref{eq3.7}. Then equation \eqref{eq3.7} is $0=1$, which is a contradiction.

If $d=0$, then $f_x^{(0)}=y_1$, which is a contradiction. Therefore, $y_1\notin \operatorname{Im}D$. Hence, the conclusion follows from the definition of Mathieu-Zhao subspace.
\end{proof}

\begin{cor}\label{cor3.3}
Let $D=\partial_x+\sum_{j=1}^n(a_j(x)y_j+b_j(x))\partial_j$ be a derivation of $K[x,y_1,\allowbreak\ldots,y_n]$. If $a_1(x),\ldots,a_n(x)$ are linearly independent and there exists $i_0\in\{1,2,\allowbreak\ldots,n\}$ such that $\deg a_{i_0}(x)\geq 1$, then $\operatorname{Im}D$ is not a Mathieu-Zhao subspace.
\end{cor}
\begin{proof}
Since $a_1(x),\ldots,a_n(x)$ are linearly independent, the equation $\sum_{i=1}^n\gamma_ia_i(x)\allowbreak=0$ has no non-zero solutions in $\mathbb{N}^n$ for $\gamma=(\gamma_1,\gamma_2,\ldots,\gamma_n)$. Thus, the conclusion follows from Theorem \ref{thm3.2}.
\end{proof}

\begin{rem}\label{rem3.4}
If $n\geq 2$, then there exists $i_0\in\{1,2,\ldots,n\}$ such that $\deg a_{i_0}(x)\geq 1$ in the case that $a_1(x),\ldots,a_n(x)$ are linearly independent over $K$. Thus, we can remove the condition that there exists $i_0\in\{1,2,\ldots,n\}$ such that $\deg a_{i_0}(x)\geq 1$ in Corollary \ref{cor3.3} if $n\geq 2$.
\end{rem}

\begin{cor}\label{cor3.5}
Let $D=\partial_x+\sum_{j=1}^r(a(x)y_j+b_j(x))\partial_j$ be a derivation of $K[x,y_1,\allowbreak\ldots,y_r]$. Then $\operatorname{Im}D$ is a Mathieu-Zhao subspace if and only if $a(x)\in K$.
\end{cor}
\begin{proof}
$``\Leftarrow"$ It follows from Lemma \ref{lem3.1} that $D$ is locally finite. Since $1\in \operatorname{Im}D$, the conclusion follows from Proposition 1.4 in \cite{8}.

$``\Rightarrow"$ If $\deg a(x)\geq 1$, then the equation $(\sum_{j=1}^r\gamma_j)a(x)=0$ has no non-zero solutions in $\mathbb{N}^r$ for $\gamma=(\gamma_1,\gamma_2,\ldots,\gamma_r)$.
It follows from Theorem \ref{thm3.2} that $\operatorname{Im}D$ is not a Mathieu-Zhao subspace. Then the conclusion follows.
\end{proof}

\begin{cor}\label{cor3.6}
Let $D=\partial_x+\sum_{j=1}^r(a(x)y_j+b_j(x))\partial_j$ be a derivation of $K[x,y_1,\allowbreak\ldots,y_r]$. Then $\operatorname{Im}D$ is a Mathieu-Zhao subspace if and only if $D$ is locally finite.
\end{cor}
\begin{proof}
The conclusion follows from Lemma \ref{lem3.1} and Corollary \ref{cor3.5}.
\end{proof}

{\bf{Acknowledgement}}:  The author is very grateful to professor Wenhua Zhao for personal communications about the Mathieu-Zhao spaces. She is also grateful to the Department of Mathematics of Illinois State University, where this paper was finished, for hospitality during her stay as a visiting scholar. The author is very grateful to the referee for some useful comments and suggestions.

\end{document}